\documentclass[11pt]{amsart}
\usepackage{amsrefs}

\def \rr {\mathbb{R}}
\def \rn {\mathbb{R}^n}

\def \sp {\mathbb{S}^p}
\def \sq {\mathbb{S}^q}

\def \crit {2^\star}
\def \ue {u_\eps}
\def \he {h_\eps}

\def \eps {\varepsilon}

\newtheorem{theorem}{Theorem}[section]
\newtheorem{proposition}[theorem]{Proposition}

\date{January 14th, 2013}

\title[Non-isolated blow-up]{Examples of non-isolated blow-up for perturbations of the scalar curvature equation on non locally conformally flat manifolds}

\author{Fr\'ed\'eric Robert}
\address{Fr\'ed\'eric Robert, Institut \'Elie Cartan, Universit\'e de Lorraine, BP 70239, F-54506 Vand{\oe}uvre-l\`es-Nancy, France}
\email{frederic.robert@univ-lorraine.fr}

\author{J\'er\^ome V\'etois}
\address{J\'er\^ome V\'etois, Universit\'e de Nice Sophia Antipolis, CNRS, LJAD, UMR 7351, F-06108 Nice, France}
\email{vetois@unice.fr}



\thanks{The authors are partially supported by the ANR grant ANR-08-BLAN-0335-01.}

\begin{document}
\begin{abstract}

Solutions to scalar curvature equations have the property that all possible
blow-up points are isolated, at least in low dimensions. This property is
commonly used as the first step in the proofs of compactness. We show that this result becomes false for
some arbitrarily small, smooth perturbations of the potential.
\end{abstract}
\maketitle
\section{Introduction and statement of the results}
Let $(M,g)$ be a smooth compact Riemannian manifold of dimension $n\geq 3$. Given a sequence $(\he)_{\eps>0}\in C^\infty(M)$, we are interested in the existence of multi peaks positive solutions $(\ue)_{\eps>0}\in C^\infty(M)$ to the family of critical equations
\begin{equation}\label{eq:1}
\Delta_g \ue+\he\ue=\ue^{\crit-1}\hbox{ in }M\hbox{ for all }\eps>0,
\end{equation}
where $\Delta_g:=-\hbox{div}_g(\nabla)$ is the Laplace-Beltrami operator, and $\crit:=\frac{2n}{n-2}$ is the critical 
Sobolev exponent. We say that the family $(\ue)_\eps$ blows up as $\eps\to 0$ if $\lim_{\eps\to 0}\Vert \ue\Vert_\infty=+\infty$. Blowing-up families to equations like \eqref{eq:1} are described precisely by Struwe \cite{struwe84} in the energy space $H_1^2$(M) : namely, if the Dirichlet energy of $\ue$ is uniformly bounded with respect to $\eps$, then there exists $u_0\in C^\infty(M)$, there exists $k\in\mathbb{N}$, there exists $k$ families $(\xi_{i,\eps})_\eps\in M$ and $(\mu_{i,\eps})_\eps\in (0,+\infty)$ such that
\begin{equation}\label{dec:struwe}
\ue=u_0+\sum_{i=1}^k\left(\frac{\sqrt{n(n-2)}\mu_{i,\eps}}{\mu_{i,\eps}^2+d_g(\cdot, \xi_{i,\eps})^2}\right)^{\frac{n-2}{2}}+o(1),
\end{equation}
where $\lim_{\eps\to 0}o(1)=0$ in $H_1^2(M)$ and $\lim_{\eps\to 0}\mu_{i,\eps}=0$ for all $i=1,...,k$. In this situation, we say that $\ue$ develops $k$ peaks when $\eps\to 0$.

\medskip\noindent We say that $\xi_0\in M$ is a blow-up point for $(\ue)_\eps$ if $\lim_{\eps\to 0}\max_{B_r(\xi_0)}\ue=+\infty$ for all $r>0$. It follows from elliptic theory that the blow-up points of a family of solutions $(\ue)_\eps$ to \eqref{eq:1} satisfying \eqref{dec:struwe} is exactly $\{\lim_{\eps\to 0}\xi_{i,\eps}/\, i=1,..,k\}$.

\medskip\noindent Following the terminology introduced by Schoen \cite{schoen1988}, $\xi_0\in M$ is an isolated point of blow-up for $(\ue)_\eps$ if there exists $(\xi_\eps)_\eps\in M$ such that
\begin{itemize}
\item $\xi_\eps$ is a local maximum point of $\ue$ for all $\eps>0$,
\item $\lim_{\eps\to 0}\xi_\eps=\xi_0$,
\item there exist $C,\bar{r}>0$ s.t. $d_g(x,\xi_\eps)^{\frac{n-2}{2}}\ue(x)\leq C$ for all $x\in B_{\bar{r}}(\xi_0)$,
\item $\lim_{\eps\to 0}\max_{B_r(\xi_0)}\ue=+\infty$ for all $r>0$.
\end{itemize}
The notion has proved to be very useful in the analysis of critical equations. 
Let $c_n:=\frac{n-2}{4(n-1)}$ and $R_g$ be the scalar curvature of $(M,g)$. Compactness for the Yamabe equation
\begin{equation}\label{eq:1Yam}
\Delta_g u + c_nR_g u= u^{2^\star-1}
\end{equation}
when $n \le 24$ (the full result is due to Kuhri--Marques--Schoen~\cite{KhuMarSch}) 
 is established by proving first that the sole possible blow-up points for \eqref{eq:1Yam} are isolated, 
see  Schoen~\cites{schoen1988,Sch3}, Li--Zhu~\cite{LiZhu}, Druet~\cite{druet:imrn}, Marques~\cite{Mar}, Li--Zhang (Theorem 1.1 in \cite{LiZha3}), and Kuhri--Marques--Schoen~\cite{KhuMarSch}). When $n\geq 25$, there are examples of non-compactness of equation \eqref{eq:1Yam} (Brendle \cite{Bre} and Brendle--Marques~\cite{BreMar}).

\medskip\noindent In this note, we address the questions to know whether or not blow-up solutions for \eqref{eq:1} do exist, and 
whether or not they necessarily have isolated blow-up points. When $\he\leq c_n R_g$, blow-up does not occur for $n\leq 5$ as shown by Druet \cite{druet:imrn} (except for the conformal class of the round sphere).  When the potential is allowed to be above the scalar curvature, blow-up is possible: we refer to Druet--Hebey \cite{druethebeyTAMS} for examples of non-isolated blow-up on the sphere with $C^1-$perturbations of the scalar curvature 
term in \eqref{eq:1Yam}, and to Esposito--Pistoia--V\'etois \cite{EspPisVet} for examples of isolated blow-up on general compact manifolds with arbitrary smooth perturbations of the scalar curvature. 
We present in this note examples of non-isolated blow-up points for smooth perturbations of the scalar curvature term in \eqref{eq:1Yam}. 
This is the subject of the following theorem.

\begin{theorem}\label{th:proof} Let $(M,g)$ be a non-locally-conformally flat compact Riemannian manifold of dimension $n\geq 6$ with positive Yamabe invariant. We fix $\xi_0\in M$ such that the Weyl tensor at $\xi_0$ is such that $\hbox{Weyl}_g(\xi_0)\neq 0$.
We let $k\geq 1$ and $r\geq 0$ be two integers. Then there exists $(h_\eps)_{\eps>0}\in C^\infty(M)$ such that $\lim_{\eps\to 0}h_\eps=c_n R_g$ in $C^r(M)$, and there exists $(\ue)_{\eps>0}\in C^\infty(M)$ a family of solutions to
$$\Delta_g\ue+h_\eps\ue=\ue^{\crit-1}\hbox{ in }M\hbox{ for all }\eps>0,$$
such that $(\ue)_\eps$ develops $k$ peaks at the blow-up point $\xi_0$. Moreover, $\xi_0$ is an isolated blow-up point if and only if $k=1$.\end{theorem}

\noindent In particular, Theorem \ref{th:proof} applies for $M:=\sp\times\sq$ ($p,q\geq 3$) endowed with the product metric. In this case, any point can be a blow-up point since the Weyl tensor never vanishes on $\sp\times\sq$.


\noindent As a consequence, when dealing with general perturbed equations like \eqref{eq:1}, one has to 
deal with the delicate situation of the accumulation of peaks at a single point. The $C^0$-theory by Druet--Hebey--Robert \cite{dhr} 
addresses this question 
in the a priori setting and $L^\infty$-norm. We refer also to Druet \cite{druet:jdg} and Druet--Hebey \cite{dh:analysis:pde} where 
the analysis of the radii of interaction of multi peaks solutions is performed.

\medskip\noindent The choice of this note is to perturb the potential $c_nR_g$ of the equation. Alternatively, one can fix the potential $c_n R_g$ and multiply
the nonlinearity $u^{2^\star-1}$ by smooth functions then leading to consider Kazdan-Warner type equations: in this slightly different context, Chen--Lin \cite{chenlinCPDE} and Brendle (private communication) have constructed non-isolated local blow-up respectively in the flat case and in the Riemannian case. 

\medskip\noindent{\bf Acknowledgements:}\,the authors express their deep thanks to\,E.\,Hebey for stimulating discussions and constant support for this project. The first author thanks C.-S. Lin for stimulating discussions and S. Brendle for communicating his unpublished result.

\section{Proofs}
The proof of Theorem \ref{th:proof} relies on a Lyapunov-Schmidt reduction. We fix $\xi_0\in M$ such that $\hbox{Weyl}_g(\xi_0)\neq 0$. It follows from the classical conformal normal coordinates theorem 
of Lee--Parker \cite{leeparker} that there exists $\Lambda\in C^\infty(M\times M)$ such that for any $\xi\in M$,
$$R_{g_\xi}(\xi)=0,\; \nabla R_{g_\xi}(\xi)=0,\text{ and }\Delta_{g_\xi} R_{g_\xi}(\xi)=\frac{1}{6}\left|\hbox{Weyl}_g(\xi)\right|^2_g,$$
where $\Lambda_\xi:=\Lambda(\xi,\cdot)$ and $g_\xi:=\Lambda_\xi^{4/(n-2)}g$. Without loss of generality, up to a conformal change of metric, we assume that $g_{\xi_0}=g$. We let $r_0>0$ be such that $r_0<i_{g_\xi}(M)$ for all $\xi\in M$ compact, where $i_{g_\xi}(M)$ is the injectivity radius of $M$ with respect to the metric $g_\xi$. We let $\chi\in C^\infty(\rr)$ be such that $\chi(t)=1$ for $t\leq r_0/2$ and $\chi(t)=0$ for $t\geq r_0$. We define a bubble centered at $\xi$ with parameter $\delta$ as:
$$W_{\delta,\xi}:=\chi(d_g(\cdot,\xi))\Lambda_\xi\bigg(\frac{\sqrt{n(n-2)}\delta}{\delta^2+d_{g_\xi}(\cdot,\xi)^2}\bigg)^{\frac{n-2}{2}}.$$
We fix an integer $k\geq 1$. Given $\alpha>1$ and $K>0$, we define the set 
\begin{multline*}
{\mathcal D}_{\alpha, K}^{(k)}(\delta):=\bigg\{((\delta_i)_{i},(\xi_i)_{i})\in (0,\delta)^k\times M^k\,\\
\left/\frac{1}{\alpha}<\frac{\delta_i}{\delta_j}<\alpha\;;\;\frac{d_g(\xi_i,\xi_j)^2}{\delta_i\delta_j}>K\,\hbox{for}\,i\neq j\right.\bigg\}.
\end{multline*}
For any $h\in C^0(M)$, we define the functional:
$$J_h(u):=\frac{1}{2}\int_M(|\nabla u|_g^2+hu^2)\, dv_g-\frac{1}{\crit}\int_M u_+^{\crit}\, dv_g$$
for all $u\in H_1^2(M)$. For $((\delta_i)_i,(\xi_i)_i)\in {\mathcal D}_{\alpha, K}^{(k)}$, we define the error
$$\hbox{$R_{(\delta_i)_i,(\xi_i)}:=\big\Vert(\Delta_g+h)\big(\sum_{i=1}^kW_{\delta_i,\xi_i}\big)-\big(\sum_{i=1}^kW_{\delta_i,\xi_i}\big)^{\crit-1}\big\Vert_{\frac{2n}{n+2}}$}$$
The classical Lyapunov-Schmidt finite-dimensional reduction yields the following:
\begin{proposition}\label{prop:gene} We fix $\alpha>1$, $\eta>0$, $C_0>0$ such that $\Vert h\Vert_\infty\leq C_0$ and $\lambda_1(\Delta_g+h)\geq C_0^{-1}$. Then there exists $K_0=K_0((M,g), \alpha, C_0, \eta)>0$, $\delta_0=\delta_0((M,g), \alpha, C_0, \eta)>0$ and $\phi\in C^1({\mathcal D}_{\alpha, K_0}^{(k)}(\delta_0) , H_1^2(M))$ such that 
\begin{itemize}
\item $R_{(\delta_i)_i,(\xi_i)}<\eta$ for all $(\delta_i)_i,(\xi_i)\in {\mathcal D}_{\alpha, K_0}^{(k)}(\delta_0)$,
\item $u((\delta_i)_i,(\xi_i)):=\sum_{i=1}^kW_{\delta_i,\xi_i}+\phi((\delta_i)_i,(\xi_i))$ is a critical point of $J_h$ iff $((\delta_i)_i,(\xi_i))$ is a critical point of$\,((\delta_i)_i,(\xi_i))\mapsto J_h(u((\delta_i)_i,(\xi_i)))\,$in ${\mathcal D}_{\alpha, K_0}^{(k)}(\delta_0)$,
\item $\Vert  \phi((\delta_i)_i,(\xi_i))\Vert_{H_1^2}=O(R_{(\delta_i)_i,(\xi_i)})$,
\item $J_h(u((\delta_i)_i,(\xi_i)_i))=J_h(\sum_{i=1}^kW_{\delta_i,\xi_i})+O(R_{(\delta_i)_i,(\xi_i)}^2)$.
\end{itemize}
Here, $|O(1)|\leq C((M,g), \alpha, C_0)$ uniformly in $ {\mathcal D}_{\alpha, K_0}^{(k)}(\delta_0)$.
\end{proposition}
\noindent This result is essentially contained in the existing litterature. It is a particular case of the general reduction theorem in  Robert--V\'etois \cite{robertvetoisPROC}. We also refer to Esposito--Pistoia--V\'etois \cite{EspPisVet} and to the general framework by Ambrosetti--Badiale \cite{AmBa} for nondegenerate critical manifolds.

\medskip\noindent From now on, we fix $((\delta_i)_i,(\xi_i)_i)\in {\mathcal D}_{\alpha, K_0}^{(k)}(\delta_0)$. Standard computations yield
\begin{multline*}
J_h\bigg(\sum_{i=1}^kW_{\delta_i,\xi_i}\bigg)=\sum_{i=1}^kJ_h(W_{\delta_i,\xi_i})+\bigg(\sum_{i\neq j}\int_M(\nabla W_{\delta_i,\xi_i},\nabla W_{\delta_j,\xi_j})_g\\
+h W_{\delta_i,\xi_i}W_{\delta_j,\xi_j}\, dv_g\bigg)-\frac{1}{\crit}\int_M\bigg(\bigg(\sum_{i=1}^k W_{\delta_i,\xi_i}\bigg)^{\crit}-\sum_{i=1}^kW_{\delta_i,\xi_i}^{\crit}\bigg)\, dv_g
\end{multline*}
and
\begin{align*}
&\int_M\bigg(\bigg(\sum_{i=1}^k W_{\delta_i,\xi_i}\bigg)^{\crit}-\sum_{i=1}^kW_{\delta_i,\xi_i}^{\crit}\bigg)\, dv_g\\
&\qquad=O\bigg(\sum_{i\neq j}\int_{W_{\delta_i,\xi_i}\leq W_{\delta_j,\xi_j}}W_{\delta_i,\xi_i}W_{\delta_j,\xi_j}^{\crit-1}\, dv_g\bigg).
\end{align*}
Choosing $K_0$ larger if necessary, there exists $c_1=c_1(\alpha,K_0)>0$ such that for any $i\neq j$ and $x\in M$ such that $W_{\delta_i,\xi_i}(x)\leq W_{\delta_j,\xi_j}(x)$, we have that $d_{g_{\xi_i}}(x,\xi_i)\geq c_1(d_g(\xi_i,\xi_j)+d_g(x,\xi_j))$. Therefore, we get that $W_{\delta_i,\xi_i}(x)\leq c_2\delta_j^{(n-2)/2}d_g(\xi_i,\xi_j)^{2-n}$ for all such $x$, for some constant $c_2=c_2(\alpha,K_0)>0$. Consequently, a rough upper bound yields
$$J_h\Big(\sum_{i=1}^kW_{\delta_i,\xi_i}\Big)= \sum_{i=1}^kJ_h(W_{\delta_i,\xi_i})+O\Big(\sum_{i\neq j}\Big(\frac{\delta_i\delta_j}{d_g(\xi_i,\xi_j)^2}\Big)^{\frac{n-2}{2}}\Big)$$
and
$$R_{(\delta_i)_i,(\xi_i)}\leq \sum_{i=1}^k\Vert (\Delta_g+h)W_{\delta_i,\xi_i}-W_{\delta_i,\xi_i}^{\crit-1}\Vert_{\frac{2n}{n+2}}+O\Big(\sum_{i\neq j}\Big(\frac{\delta_i\delta_j}{d_g(\xi_i,\xi_j)^2}\Big)^{\frac{n-2}{4}}\Big)$$
uniformly in ${\mathcal D}_{\alpha, K_0}^{(k)}(\delta_0)$. Moreover, see Proposition 2.3 in Esposito--Pistoia--V\'etois \cite{EspPisVet}, we have that
\begin{multline*}
J_h(W_{\delta,\xi})=\frac{K_n^{-n}}{n}\bigg(1+\frac{2(n-1)}{(n-2)(n-4)}(h-c_nR_g)(\xi)\delta^2\\
+O(\Vert h-c_n R_g\Vert_{C^{1}})\delta^{3}\\
-|Weyl_g(\xi)|_g^2\left\{\begin{array}{ll}
\frac{1}{64}\delta^4\ln\frac{1}{\delta} +O(\delta^4)& \hbox{ when }n=6\\
\frac{1}{24(n-4)(n-6)}\delta^4+O(\delta^5) & \hbox{ when }n\geq 7
\end{array}\right.\bigg)
\end{multline*}
and
\begin{multline*}
\Vert (\Delta_g+h)W_{\delta,\xi}-W_{\delta,\xi}^{\crit-1}\Vert_{\frac{2n}{n+2}}\\
\leq C\delta^2\left\{\begin{array}{ll}
1+\Vert h-c_n R_g\Vert_{C^{0}}\left(\ln\frac{1}{\delta}\right)^{2/3}& \hbox{ when }n=6\\
\sqrt{\delta}+\Vert h-c_n R_g\Vert_{C^{0}}& \hbox{ when }n\geq 7.
\end{array}\right.
\end{multline*}
Here again, $|O(1)|\leq C((M,g), \alpha, C_0)$ uniformly in ${\mathcal D}_{\alpha, K_0}^{(k)}(\delta_0)$.

\medskip\noindent We now choose the $(\delta_i),(\xi_i)'s$ and the function $h$. For any $\eps>0$, we let $\delta_\eps>0$ be such that
$$\delta_\eps^2\ln\frac{1}{\delta_\eps}=\eps\hbox{ when }n=6\hbox{ and }\delta_\eps^2=\eps\hbox{ when }n\geq 7.$$
We let $H\in C^\infty(\rn)$ be such that 
\begin{itemize}
\item $H(x)=-1$ for all $|x|>2$,
\item $H$ admits $k$ distinct strict local maxima at $p_{i,0}\in B_1(0)$ for $i=1,...,k$,
\item $H(p_{i,0})>0$ for all $i=1,...,k$.
\end{itemize}
We let $\tilde{r}>0$ be such that for any $i\in \{1,...,k\}$, the maximum of $H$ on $B_{2\tilde{r}}(p_{i,0})$ is achieved exactly at $p_{i,0}$ and such that $|p_{i,0}-p_{j,0}|\geq 3\tilde{r}$ for all $i\neq j$. We let $(\mu_\eps)_\eps\in (0,+\infty)$ be such that $\lim_{\eps\to 0}\mu_\eps=0$ and
$$\left(|\ln\eps|\right)^{-1/4}=o(\mu_\eps)\hbox{ when }n=6\hbox{ and }\eps^{\frac{n-6}{2(n-2)}}=o(\mu_\eps)\hbox{ when }n\geq 7,$$
where both limits are taken when $\eps\to 0$. As one can check, $\delta_\eps=o(\mu_\eps)$ when $\eps\to 0$. We define 
$$h_\eps(x):=c_n R_g(x)+\eps H\left(\mu_\eps^{-1}\hbox{exp}_{\xi_0}^{-1}(x)\right)\hbox{ for all }x\in M.$$
Here, the exponential map is taken with respect to the metric $g$ and after assimilation to $\rn$ of the tangent space at $\xi_0$: this definition makes sense for $\eps>0$ small enough. For $(t_i)_i\in (0,+\infty)^k$ and $(p_i)_i\in (\rn)^k$, we define
$$\tilde{u}_\eps((t_i)_i,(p_i)_i):=u\left( (t_i \delta_\eps)_i, (\hbox{exp}_{\xi_0}(\mu_\eps p_i))_i\right)\hbox{ with }h\equiv h_\eps.$$
The above estimates and the choice of the parameters yield
\begin{equation}\label{lim:J}
\lim_{\eps\to 0}\frac{J_{h_\eps}(\tilde{u}_\eps((t_i)_i,(p_i)_i))-k\frac{K_n^{-n}}{n}}{\eps\delta_\eps^2}=\sum_{i=1}^kF_n(t_i, p_i)
\end{equation}
in $C^0_{loc}((0,+\infty)^k\times \prod_{i=1}^k B_{r}(p_{i,0}))$, where  
$$F_n(t,p):=\frac{2(n-1)}{(n-2)(n-4)}H(p)t^2-d_n |Weyl_g(\xi_0)|_g^2t^4$$
for $(t,p)\in (0,+\infty)\times \rn$, with $d_6=\frac{1}{64}$ and $d_n:=\frac{1}{24(n-4)(n-6)}$ for $n\geq 7$. As easily checked, up to choosing the $t_i's$ in suitable compact intervals $I_1,...,I_k$, the right-hand-side of \eqref{lim:J} has a unique maximum point in the interior of  $ \prod_{i=1}^k I_i\times \prod_{i=1}^k B_{\tilde{r}}(p_{i,0})$. As a consequence, for $\eps>0$ small enough, $J_{h_\eps}(\tilde{u}_\eps((t_i)_i,(p_i)_i))$ admits a critical point, $((t_{i,\eps})_i,(p_{i,\eps})_i)\in (\alpha,\beta)^k\times \prod_{i=1}^k B_{\tilde{r}}(p_{i,0})$ for some $0<\alpha<\beta$ independent of $\eps$. Defining $\xi_{i,\eps}:=\hbox{exp}_{\xi_0}(\mu_\eps p_{i,\eps})$ for all $i=1,...,k$, there exists $c_0>0$ such that $d(\xi_{i,\eps}, \xi_{i,\eps})\geq c_0\mu_\eps$ for all $i\neq j\in\{1,...,k\}$ and all $\eps>0$ small enough. Defining $u_\eps:=\tilde{u}_\eps((t_{i,\eps})_i,(p_{i,\eps})_i)$, it follows from Proposition \ref{prop:gene} and the strong maximum principle that
$$\Delta_g u_\eps+h_\eps u_\eps=u_\eps^{\crit-1}\hbox{ in }M$$
for $\eps>0$ small enough. In addition to the hypotheses above, we require that $\eps=o(\mu_\eps^r)$ when $\eps\to 0$, which yields $\lim_{\eps\to 0}h_\eps=c_n R_g$ in $C^r(M)$.

\medskip\noindent We prove that $(u_\eps)_\eps$ develops no isolated blow-up point when $k\geq 2$. We argue by contradiction. Moser's iterative scheme yields the convergence of $\ue$ to $0$ in $C^2_{loc}(M\setminus\{\xi_0\})$. We then get that the isolated blow-up point is $\xi_0$, and thus that there exists $r_1>0$ and $(\xi_\eps)_\eps\in M$ such that $\lim_{\eps\to 0}\xi_\eps=\xi_0$ and there exists $C>0$ such that
\begin{equation}\label{ineq:simple}
d_g(x,\xi_\eps)^{\frac{n-2}{2}}\ue(x)\leq C\hbox{ for all }\eps>0\hbox{ and }x\in B_{r_1}(\xi_0).
\end{equation}
For any $i=1,..,k$, we recall that $\xi_{i,\eps}:=\hbox{exp}_{\xi_0}(\mu_\eps p_{i,\eps})$ and we define
$$\tilde{u}_{i,\eps}(x):=(\delta_\eps t_{i,\eps})^{\frac{n-2}{2}}u_\eps(\hbox{exp}_{\xi_{i,\eps}}(\delta_\eps t_{i,\eps} x))$$
for all $|x|<r_0/(2\delta_\eps t_{i,\eps})$. It follows from standard elliptic theory that
\begin{equation}\label{lim:U}
\lim_{\eps\to 0}\tilde{u}_{i,\eps}=\bigg(\frac{\sqrt{n(n-2)}}{1+|\cdot |^2}\bigg)^{\frac{n-2}{2}}\hbox{ in }C^2_{loc}(\rn).
\end{equation}
Moreover, if $\delta_\eps =o(d_g(\xi_{i,\eps},\xi_\eps))$ when $\eps\to 0$, inequality \eqref{ineq:simple} yields the convergence of $\tilde{u}_{i,\eps}$ to $0$ in $C^0_{loc}(\rn)$: a contradiction to \eqref{lim:U}. Therefore, $d_g(\xi_\eps, \xi_{i,\eps})=O(\delta_\eps)$ when $\eps\to 0$ for all $i=1,...,k$, and then $d_g(\xi_{i,\eps}, \xi_{j,\eps})=O(\delta_\eps)=o(\mu_\eps)$ when $\eps\to 0$ for all $i\neq j$. This contradicts the fact that $d_g(\xi_{i,\eps}, \xi_{j,\eps})\geq c_0\mu_\eps$ when $k\geq 2$. This proves the non-simpleness when $k\geq 2$.

\end{document}